\documentclass[aip,cha,preprint]{revtex4-1}
\usepackage{amsmath,amssymb,graphicx,url}

\newcommand{\RR}{\mathbb{R}}

\newcommand{\ZZ}{\mathbb{Z}}
\newcommand{\CC}{\mathbb{C}}
\newtheorem{Theorem}{Theorem}
\newtheorem{Proposition}[Theorem]{Proposition}

\begin{document}
\title{New relatives of the Sierpinski gasket} 
\author{Christoph Bandt}
\affiliation{Institute of Mathematics, University of Greifswald, 17487 Greifswald, Germany} 
\email{bandt@uni-greifswald.de}
\author{Dmitry Mekhontsev}
\affiliation{Sobolev Institute of Mathematics, 4 Acad.\  Koptyug avenue, 630090 Novosibirsk Russia}
\email{mekhontsev@gmail.com}
\date{\today}
%\subjclass[2010]{28A80}

\begin{abstract}
By slight modification of the data of the Sierpinski gasket, keeping the open set condition fulfilled, we obtain self-similar sets with very dense parts, similar to fractals in nature and in random models. This is caused by a complicated structure of the open set and is revealed only under magnification. Thus the family of self-similar sets with separation condition is much richer and has higher modelling potential than usually expected. An interactive computer search for such examples and new properties for their classification are discussed. 
\end{abstract}

\maketitle
 
\section{Introduction}\label{intro}
\begin{figure}[h!] 
\begin{center}
\includegraphics[width=0.9\textwidth]{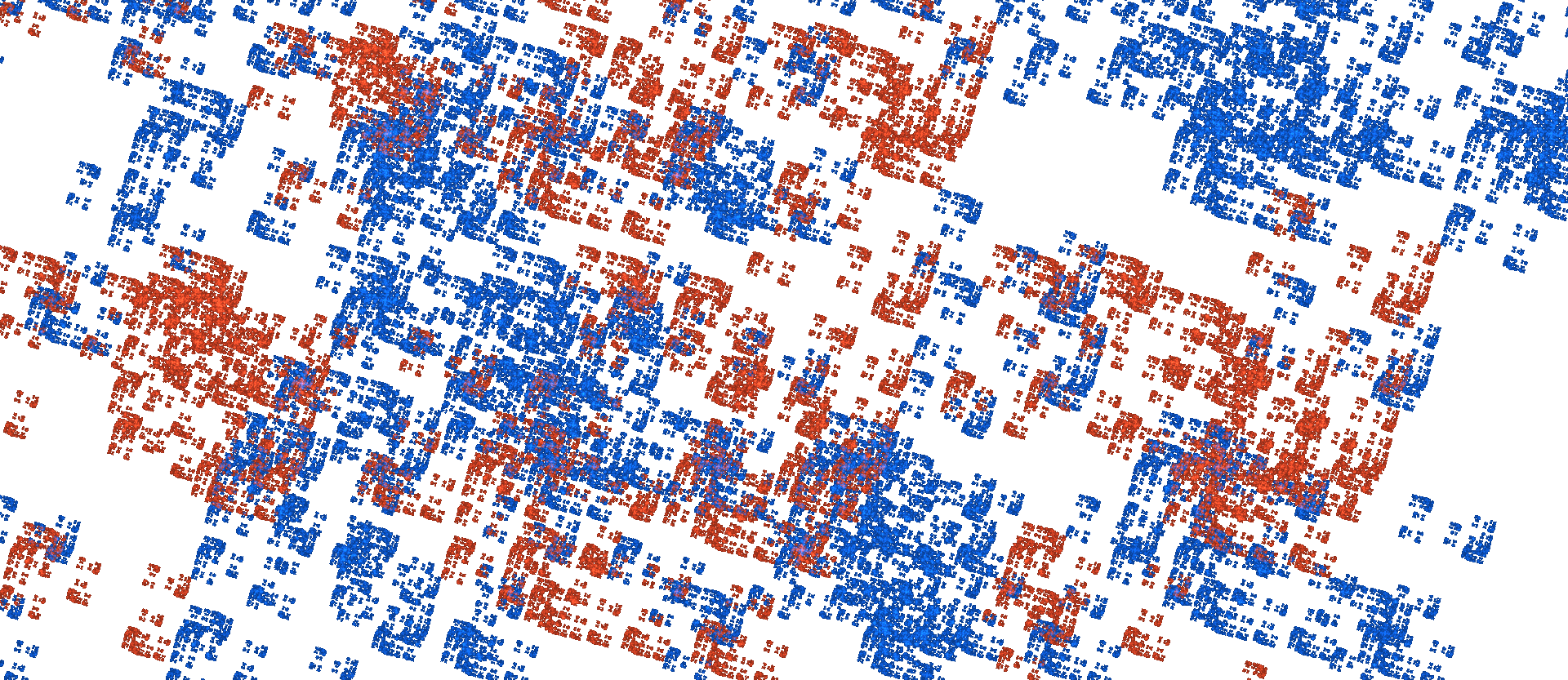}
\end{center}
\caption{``Patches''. A detail of the self-similar Cantor set in Figure \ref{examples} generated by the similarity maps $f_0(z)=\frac{iz}{2}-3-5i,\, f_1(z)=\frac{z}{2}+1+\frac{5i}{2},\, f_2(z)=\frac{-iz}{2}+1+3i$ in $\CC .$ }\label{patches}
\end{figure}

\begin{figure}[h!] 
\begin{center}
\includegraphics[width=0.9\textwidth]{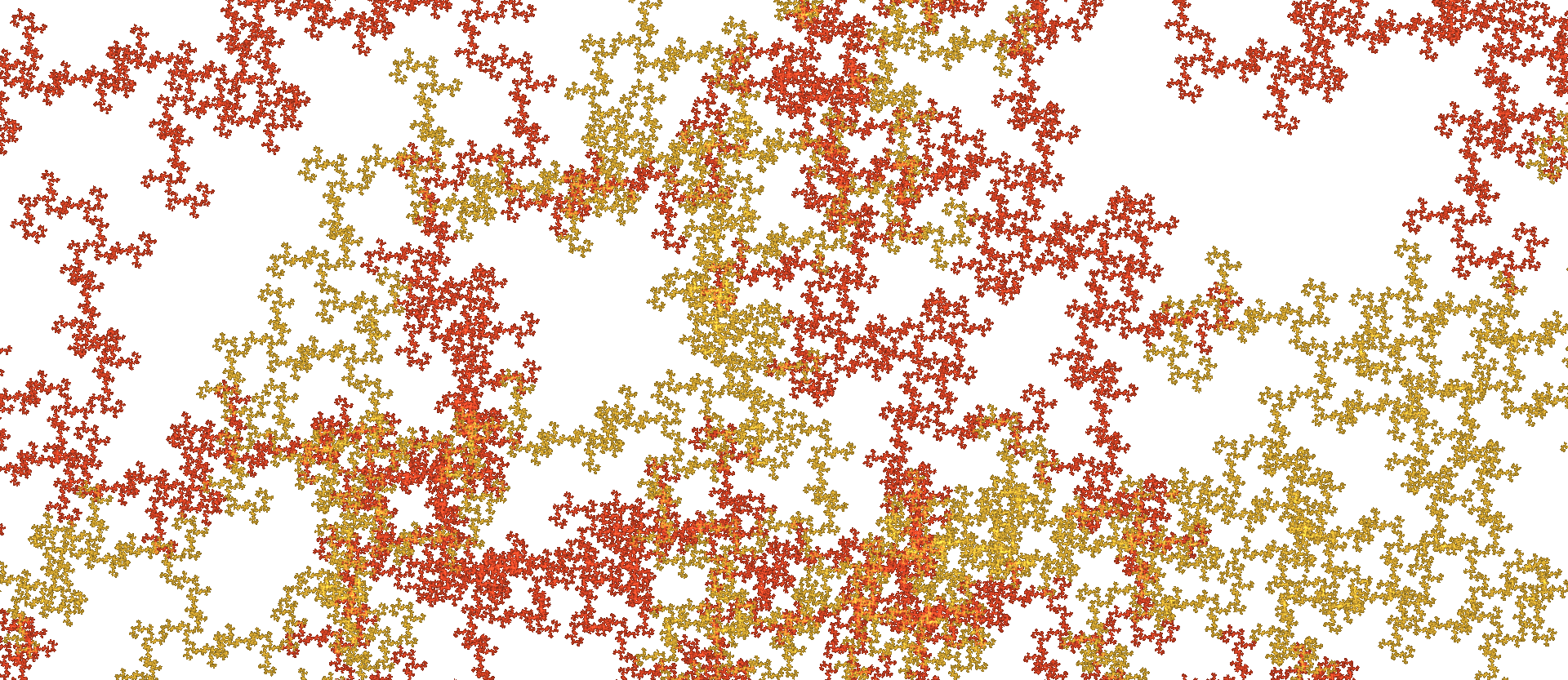}
\end{center}
\caption{Magnification of the set ``Thicket'' in Figure \ref{examples} generated by the similarity mappings $f_0(z)=\frac{iz}{2} -1-\frac{i}{2},\, f_1(z)=\frac{iz}{2}-\frac{3i}{2},\, f_2(z)=-\frac{z}{2}+\frac{i-1}{2}$.}\label{thicket}
\end{figure}

\begin{figure}[h!] 
\begin{center}
\includegraphics[width=0.9\textwidth]{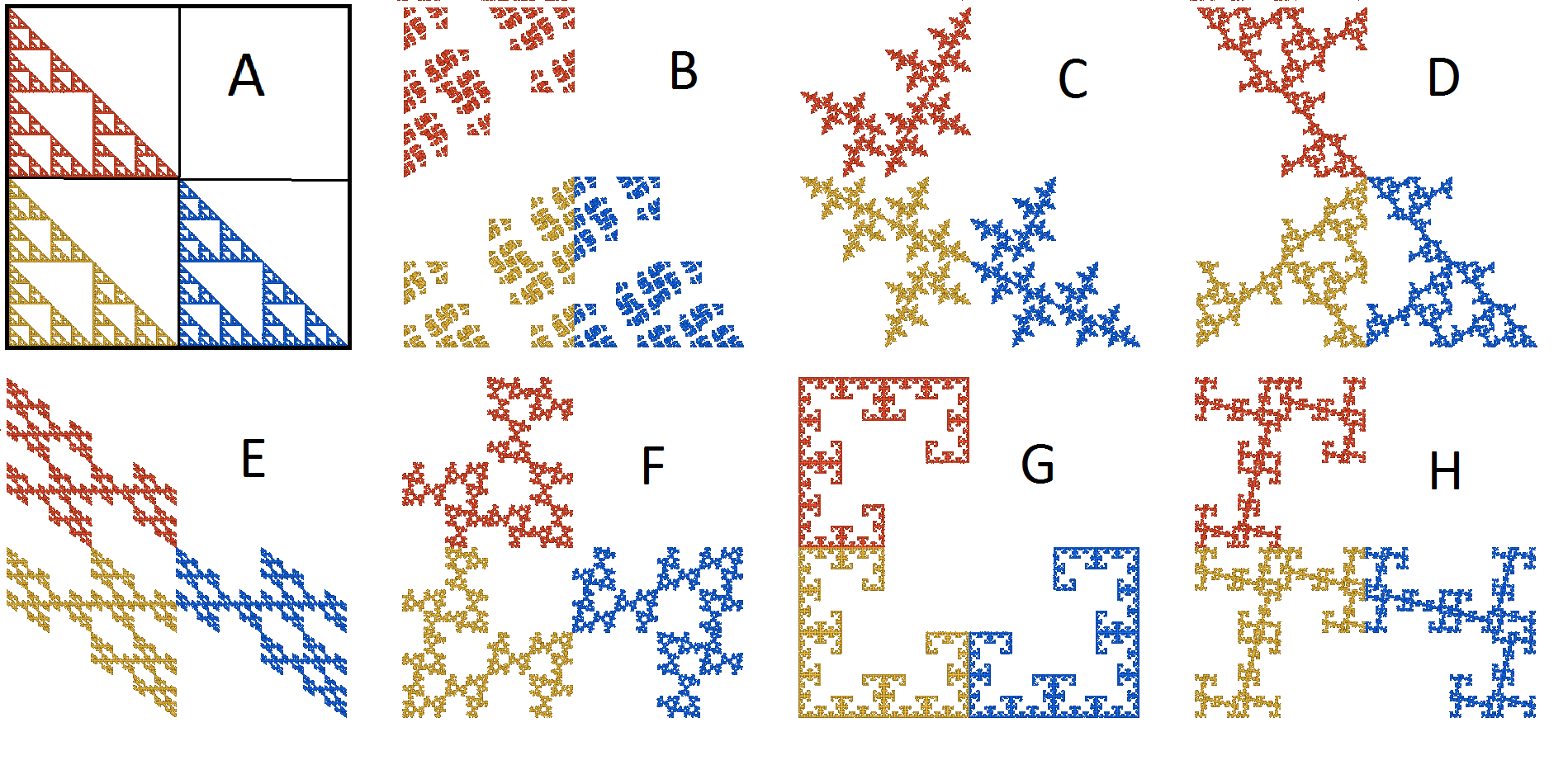}
\end{center}
\caption{Well-known fractals \cite{bar}. A. The Sierpinski gasket, B. Cantor set, C. disconnected set containing infinitely many intervals, D-H. connected examples. The square can be taken as open set $U.$ Only for the sets F and H, magnification will provide new details.}\label{fsquares}
\end{figure}

Self-similar sets with open set condition are the simplest mathematical fractals, and the Sierpinski gasket, shown in Figure \ref{fsquares}A, is their best-known prototype. It seems obvious that they do not change their structure when magnified, and that they are much more regular than fractals in nature, apart from Barnsley's fern \cite{bar} and Romanesco cauliflower. In this note we show that both these statements are wrong. Many modifications of the Sierpinski gasket, like Figures \ref{patches} and \ref{thicket}, show their structure only under strong magnification, and their complexity is much larger than expected from the definition. Below we give more examples, and discuss their computer-assisted construction and their properties.

\paragraph{Basic definitions.}
A non-empty compact set $A\subset \mathbb R^d$ is called a \emph{self-similar set} with respect to $m$ contractive similarity mappings $f_1,...,f_m$ if it fulfils the equation 
\begin{equation}\label{selsi} 
A=f_1(A)\cup ...\cup f_m(A) \, ,
\end{equation}
which means that $A$ is a union of $m$ similar copies of itself. 
A basic theorem of Hutchinson says that for given contracting maps $f_1,...,f_m,$ equation \eqref{selsi} has a unique solution $A$ \cite{bar,Fal}.  The family $\{ f_1,...,f_m\}$ was termed iterated function system (IFS) by Barnsley since the set $A$ can be obtained by successive iteration of the maps $f_k,$ see  \cite{bar,PJS}.
It turns out that only special IFS provide a nice structure of $A.$ One usually assumes the open set condition (OSC) which says that there is an open set $U$ such that
\begin{equation}\label{osc} 
f_j(U)\subset U\quad\mbox{ and }\quad f_j(U)\cap f_k(U)=\emptyset \quad\mbox{ for } j,k=1,...,m \mbox{ with }  j\not= k\, .
\end{equation}
In Figure \ref{fsquares}, $U$ can be taken as an open square which contains all points of $A$ with exception of a few boundary points.
The sets $f_j(U),$ and their smaller subsets $f_{j_1}(f_{j_2}(U)),$ and $f_{j_1}f_{j_2}\cdots f_{j_n}(U)$ provide nets of
disjoint open sets which prescribe the positions of the smaller pieces $f_{j_1}f_{j_2}\cdots f_{j_n}(A)$ of $A.$ This was the idea of P.A.P. Moran \cite{Mo} who introduced the OSC in 1946. In Figure \ref{fsquares} the smaller pieces are well recognizable since they do not overlap. This holds whenever $U$ is convex, or a polygon. In Figures \ref{patches} and \ref{thicket}, however, apparent overlap is possible since $U$ is more complicated.

\paragraph{Our class of IFS.}
The Sierpinski gasket is generated by three maps $f_k(z)=\frac12 (z+c_k)$ where the $c_k$ are non-collinear points. We consider related mappings in the complex plane $\CC .$ We always take three maps of the form
\begin{equation}\label{maps}
  \textstyle f_k(z)=\frac12 \cdot s_k(z+c_k) \quad \mbox{ for } k=0,\, 1,\, 2\,
\end{equation}
where $c_k=a_k+i\cdot b_k$ is a complex number with integers $a_k, b_k,$ and each $s_k$ is an element of 
\[ S=\{ id,\,  s(z)=i\cdot z,\, s^2(z)=-z,\,  s^3(z)=-i\cdot z\} \, .\]
Thus all maps are integer translations
combined with a rotation around $0, \pm 90$ or $180$ degrees and a homothety contracting by $\frac12 .$

\paragraph{Contents of the paper.}
In Section \ref{well} we briefly analyze eight well-known examples \cite{bar,PJS,FC} which are  obtained by taking a square as separating open set $U.$  Although they differ in geometric and topological properties, the shape of $U$ imposes severe restrictions on the structure of the fractals $A.$

For our new examples, the open set was not prescribed. Its existence is checked algebraically, with help of a computer. 
The second author's IFS tile finder package \cite{Me} was used for screening our parameter space of IFS. This package works for much more complicated IFS in two or three dimensions without much additional effort. We confine ourselves to particular simple IFS defined by 20 bytes of data, showing that even slight modifications of the well-known Sierpinski gasket provide new insights. These fractals are of course still far from modelling nature. They look a bit rectangular since only right angles were used in rotations. All have the same fractal dimension.

Six selected examples are presented in Section \ref{exa}. In Section \ref{open} we discuss their open sets $U$ and derive a new construction for $U.$ In Section \ref{neigh} we introduce our main tool, the neighbor graph, and determine various properties of our examples. The final Section \ref{last} deals with problems and options of the computer search. All  examples and the software for their analysis are available on the web site \cite{Me}.

\section{Some well-known examples} \label{well}
\paragraph{Definition.} Our class of IFS includes the examples of Figure \ref{fsquares}, where we can arrange the square so that its center is at zero and three vertices are $c_0=-1-i,\, c_1=1-i,$ and $c_2=-1+i.$ Then we can take $f_k(z)=\frac12 (s_k(z) +c_k).$ Thus to get a small copy we first rotate the big square around a multiple of $90^o,$ then translate it by $c_k$ and contract it so that it comes to one of the small subsquare with vertex $c_k.$ 

It is easy to identify the $s_k$ in each of these examples. With four rotations to choose for $k=0,1,2,$ we have 64 cases. We identify each set with its mirror image at the line $y=x.$ Then we are left with 36 sets, eight of which are symmetric, like Figure \ref{fsquares}A and G. See Barnsley \cite[Figure VII.200]{bar}. A larger family, with reflections of the square included as symmetries $s,$ was studied by Peitgen, J\"urgens and Saupe \cite{PJS} and by Falconer and O'Connor \cite{FC}. Lau, Luo and Rao considered `fractal squares' with an arbitrary choice among $n\times n$ subsquares for $n>2,$ but without rotations \cite{LLR}. Self-similar subdivisions of a triangle have also been considered.

In our pictures $A_0=f_0(A)$ is drawn in yellow, $A_1$ in blue, and $A_2$ in red color. All examples in Figure \ref{fsquares} have the same fractal dimension $d=\frac{\log 3}{\log 2}\approx 1.58 .$ Nevertheless, they differ a lot in topology, and in the structure of intersections of pieces, as Table \ref{fsq} indicates. \vspace{2ex} . 

\begin{table}[h!] 
\begin{center}
\begin{tabular}{|l|c|c|c|c|c|c|c|c|}\hline
name &\bf A&\bf B&\bf C&\bf D&\bf E&\bf F&\bf G&\bf H\\ \hline
topology&conn.&Cantor&disconn.&conn.&conn.&conn.&conn.&conn.\\ \hline
holes&yes&no&no&yes&yes&yes&no&yes\\ \hline
line segments&yes&no&yes&yes&yes&yes&yes&no\\ \hline
proper nbs.&6&11&2&10&5&9&7&23 \\ \hline
finite nbs.&6&6&2&10&5&2&3&9  \\ \hline
boundary dim&0&$\frac{d}{2}$&0&0&0&$\frac{1}{4}$&1&$\frac{1}{3},\, \frac{2}{3}$  \\ \hline
max degree&3&5&2&4&4&3&7&6  \\ \hline
neighborhoods&6&14&1&14&10&15&12&42  \\ \hline
\end{tabular}\vspace{-2ex} 
\end{center}
\caption{Properties of the fractal squares in Figure \ref{fsquares}, explained in the text.}  \vspace{-2ex}
\label{fsq}\end{table}

\paragraph{Topological properties.} 
`Topology' refers to connectedness. $A$ is connected if and only if there is a piece $A_k$ which intersects both of the other pieces. This is fulfilled for all examples except B and C. In C, only the pieces $A_0,A_1$ intersect in a single point, and this is sufficient to generate an interval, namely the self-similar set $I$ with respect to $f_0,f_1.$ 
As a consequence, there are infinitely many intervals in C obtained from $I$ by repeated application of the maps $f_k.$ (This corrects a remark in \cite[p. 251]{PJS}.) Such intervals appear whenever two of the maps have rotation angles 0 or 180 degrees. 

In example B, the intersection $A_0\cap A_1$ is a Cantor set with dimension $\frac{d}{2}\approx 0.79$, and still there are no connected subsets with more than one point in this example (cf. Section \ref{neigh}).
The maximal dimension of intersection sets of pieces will be termed boundary dimension. It is an important parameter to distinguish our examples.

A set is simply connected if it does not enclose holes. 
In example G, the intersections of neighboring pieces have the largest possible dimension 1 since they are intervals, and still the set is simply connected. In examples A, E and D, two pieces intersect in at most two points, and the sets surround infinitely many holes.

\paragraph{Boundary and neighbor structure.} 
F and H are the only examples where magnification of the overlap region, in particular $A_0\cap A_2,$ shows new details.  Our new examples will cause even more surprise when magnified. H does not contain any intervals. F contains intervals like C, but it seems not so clear how we can connect any two points by a rectifiable curve. The boundary dimension is $\frac{1}{4}$ for F, and  there are boundary sets of dimension $\frac{1}{3}$ and $\frac{2}{3}$ in example H. 

To prove these statements, we need the types of intersections of pieces of $A,$ as defined in Sections \ref{open}, \ref{neigh}. We call them `proper neighbors', and `finite neighbors' when the intersection is a finite set. The Sierpinski gasket has three vertices where a neighbor could be appended when we extend the fractal construction outwards. At the upper vertex we can append a neighbor to the left or to the right. For the other vertices we also have two cases, so altogether there are 6 finite neighbors. Neighboring pieces with infinite intersection exist in examples B, F, G, and H.

An important parameter in Table \ref{fsq} is the maximum degree which counts the maximal number of neighbors which actually occur at a single small piece of $A.$  By OSC, this number must be finite \cite{Sch}, and in our new examples it will be fairly large. Here it ranges between 0 and 7. We shall also count the number of neighborhoods, that is, combinations of neighbors which appear at small subpieces everywhere in $A.$ Most of our new examples admit thousands of neighborhoods.
%, cf. Table \ref{paras} below.  

\paragraph{Restrictions.} 
Although there is considerable geometric variety in Figure \ref{squares}, the choice of the similarity maps and the square $U$ impose restrictions:
\begin{enumerate}
\item All IFS have three similitudes with factor $\frac12$ and rotation angles of $0, 180, \pm 90$ degrees. As a consequence, the dimension is $d=\frac{\log 3}{\log 2}\approx 1.585 .$
\item The intersection of two pieces must always be a subset of an interval. 
\item Any piece of the fractal has at most 8 neighbors, and at most 4 which intersect in more than one point.
\end{enumerate}

Prescribing the square as open set $U$ is done to construct examples easily, but it is not natural. The similarity of branches of a tree is not caused by reservation of disjoint regions of equal size for different pieces. In nature they will often intermingle. 
Below we shall avoid the restrictions 2 and 3 by checking the OSC algebraically. We keep restriction 1 to keep our presentation simple.

\begin{figure}[h!] 
\begin{center}
\includegraphics[width=0.99\textwidth]{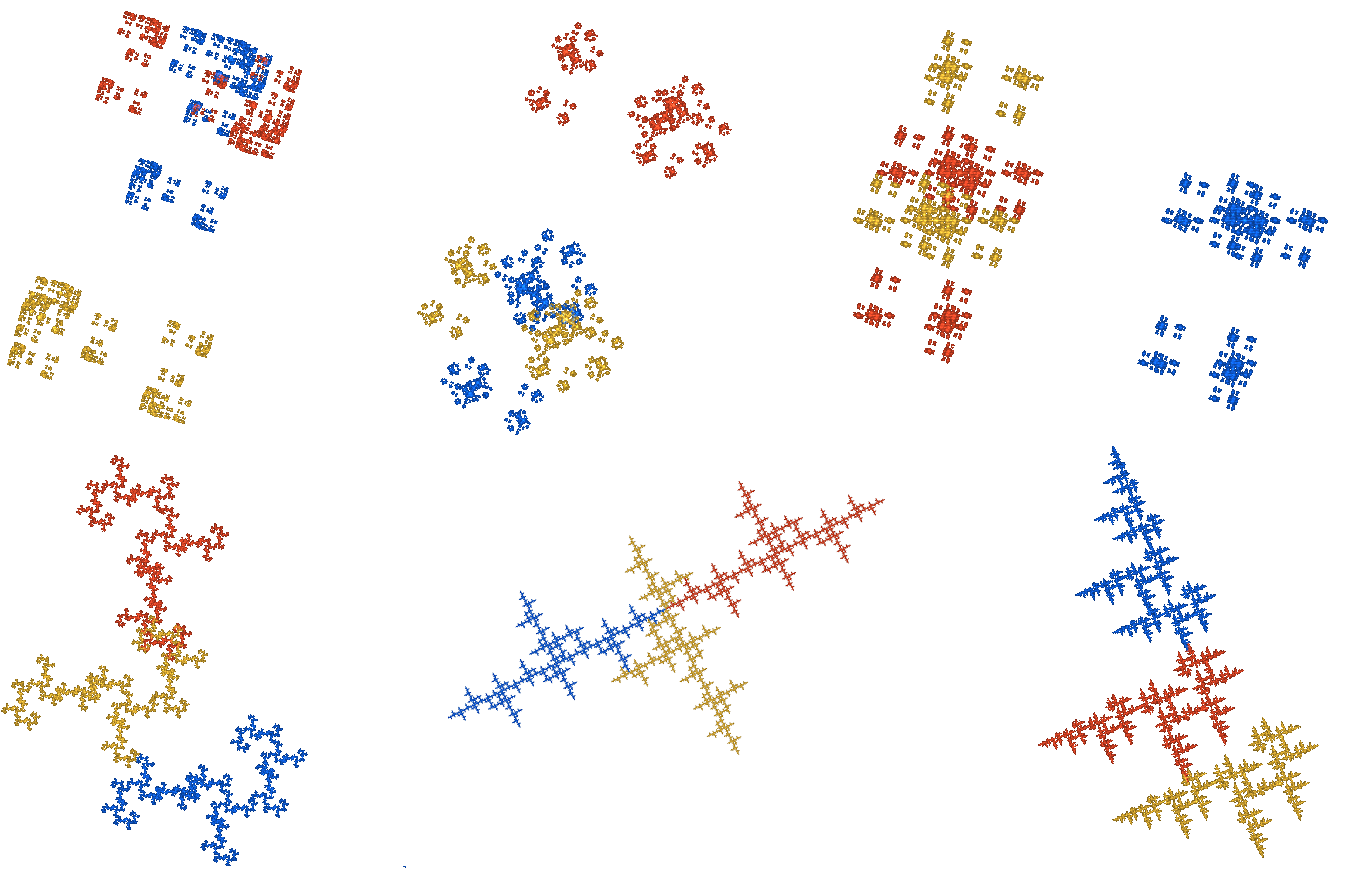}
\end{center}
\caption{The Cantor sets ``Patches'', ``Bubbles'', and ``Fireworks'' are shown in the upper row, the connected sets ``Thicket'', ``Crossings'', and ``Forest'' in the lower row. In spite of apparent overlaps, the open set condition is fulfilled. For details see Figures \ref{patches}, \ref{thicket}, \ref{firework}, and \ref{crossing}.}\label{examples}
\end{figure}  

\section{Six new examples} \label{exa}
We explore the family of maps by computer, using the program IFStile which is freely available at \cite{Me}. Roughly speaking, the program performs a kind of random walk in the parameter space $\Pi$ of all possible IFS:
\begin{equation} 
\Pi = \{ (s_0, s_1, s_2, c_0, c_1, c_2)\, |\, s_k\in S, c_k=a_k+i\cdot b_k, a_k,b_k\in \ZZ\}\, . 
\label{pi}\end{equation}
In each step, either some $s_k$ is changed, and/or some $a_k$ or $b_k$ is changed by $\pm 1.$ Many versions of the random walk are possible, see Section \ref{last}. Each resulting IFS is analyzed, and is added to the list of examples if it fulfils the open set condition and has interesting properties. Figure \ref{examples} shows a sample of six self-similar sets, selected from tens of thousands of examples obtained during two hours of computer search.

\begin{figure}[h!] 
\begin{center}
\includegraphics[width=0.9\textwidth]{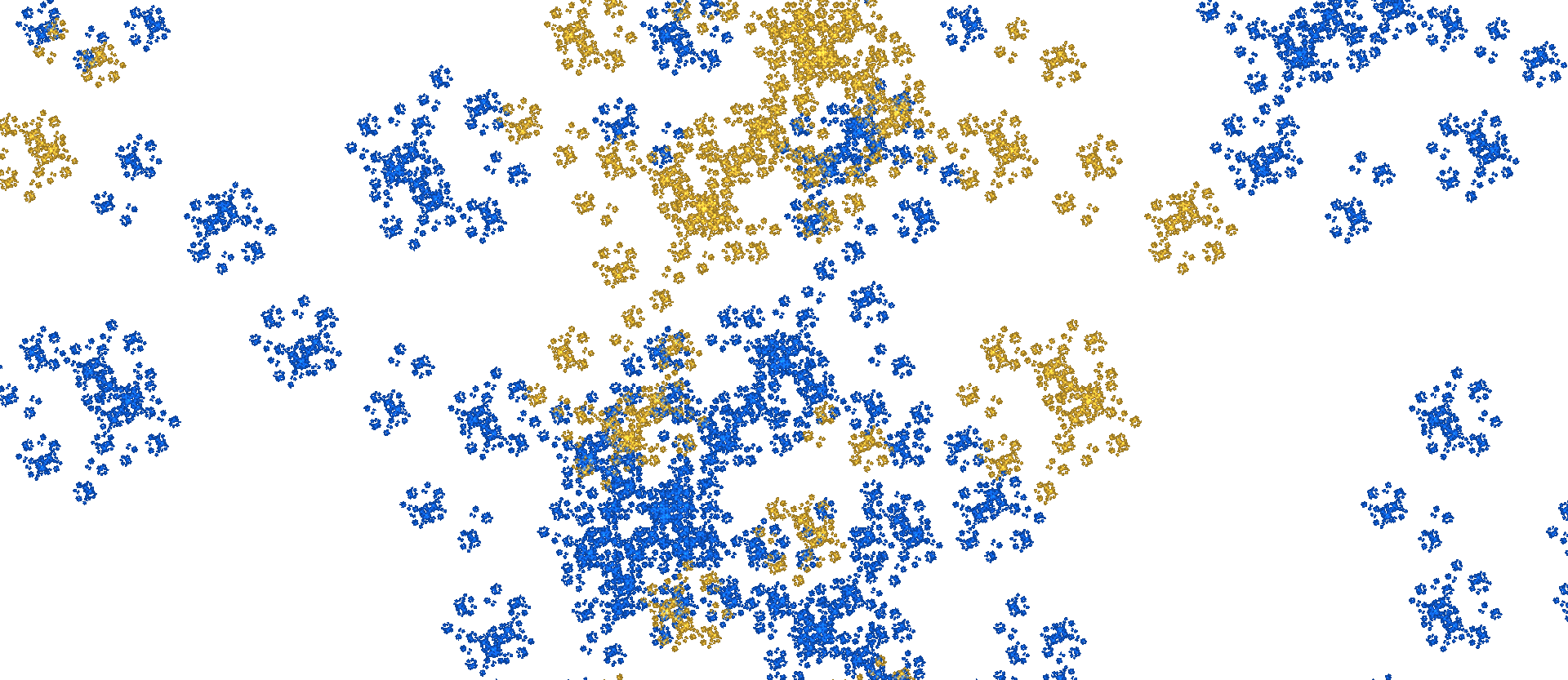}\vspace{4ex}\\
\includegraphics[width=0.9\textwidth]{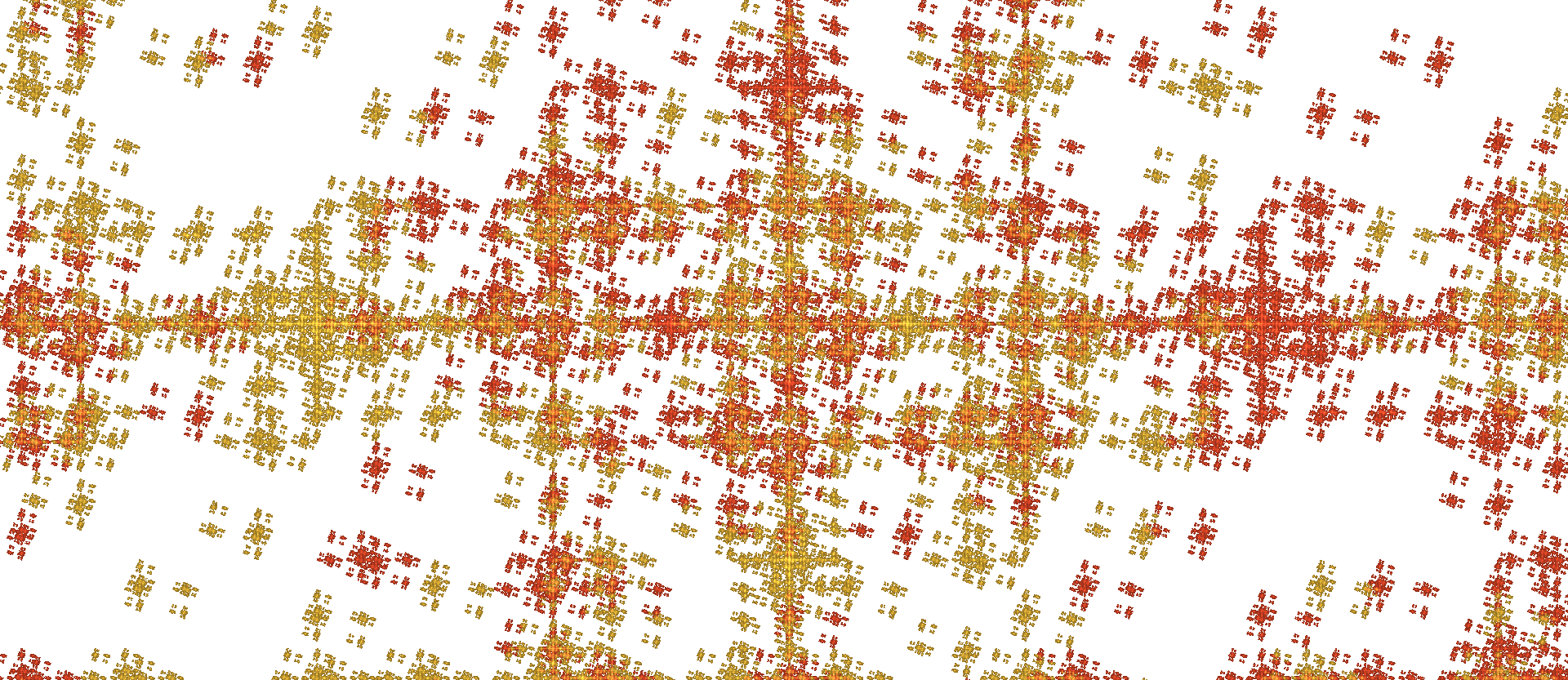}
\end{center}
\caption{Details of two Cantor sets in Figure \ref{examples}. ``Bubbles'' is generated by $f_0(z)=\frac{iz}{2} -1-i,\, f_1(z)=-\frac{z}{2}+1+\frac{3i}{2},\, f_2(z)=\frac{iz}{2}+i,$ and ``Fireworks'' by $f_0(z)=\frac{iz}{2},\, f_1(z)=\frac{-iz}{2} +3,\, f_2(z)=\frac{-iz}{2}+\frac{i}{2}.$ }\label{firework}
\end{figure}  

\paragraph{Magnification is crucial.}
It is quite natural that the structure of these objects becomes visible only under magnification. A photo of a tree on a field, seen from some distance, does not show the self-similar appearance of branches and leaves. You must go nearer to recognize  structure. Details look different at various places, but are so tightly related that we can guess the species from each part.  This also holds for abstract self-similar sets. They are defined as compact sets, but their structure is revealed in a collection of detail pictures. 

Here we focus on examples which look harmless as compact sets and become dense and intricate in magnification. However, our magnifications contain isolated copies of the original set. So by further zooming at the right places we obtain thin views again. Who observed a cloud from an airplane, or a tree in the garden, will confirm that nature mixes thin and dense structure. 

In Figure \ref{examples}, $A_0=f_0(A)$ is again drawn in yellow, $A_1$ in blue, and $A_2$ in red color.The angles of rotation of the $f_k$ can be read from the picture. The formula for the $f_k$ is given in the captions of magnified figures. Magnification was performed within the overlap region, at places where many small pieces cluster together. Below, we will argue that in this way we see structure of $A$ which regularly repeats when we continue zooming. Places with little overlap can belong to the dynamical boundary of $A$ and need not repeat in the magnification flow when we zoom around a typical point of $A.$

\paragraph{Microsets.}
There is a mathematical theory of `microsets' (Furstenberg, Hochman) or 'tangential measure distributions' (Bandt, Graf, M\"orters and Preiss) which defines the `space of all limit views of a fractal' in a compelling way. For a self-similar fractal of finite type (see Section \ref{neigh}), we need not go to the limit, and the space of views is a compact probability space. See \cite{Ba3, Ho} for details and more references. Roughly speaking, we take all magnifications which are inside the open set $U,$ but not inside the $f_k(U),$ because there everything will repeat. The resulting collection of pictures will not depend on the choice of $U.$ We take this rough idea as a guideline.

The upper row of Figure \ref{examples} contains disconnected examples: only two of the three pieces intersect. Magnification of the overlap of these two pieces in Figures \ref{patches} and \ref{firework} indicate that there could be connected subsets.  Careful further magnification verifies that no connected subsets exist: these three examples are Cantor sets. From a physicist's point of view, we are near to a phase transition with percolation -- connections are going to emerge.

\begin{figure}[h!] 
\begin{center}
\includegraphics[width=0.9\textwidth]{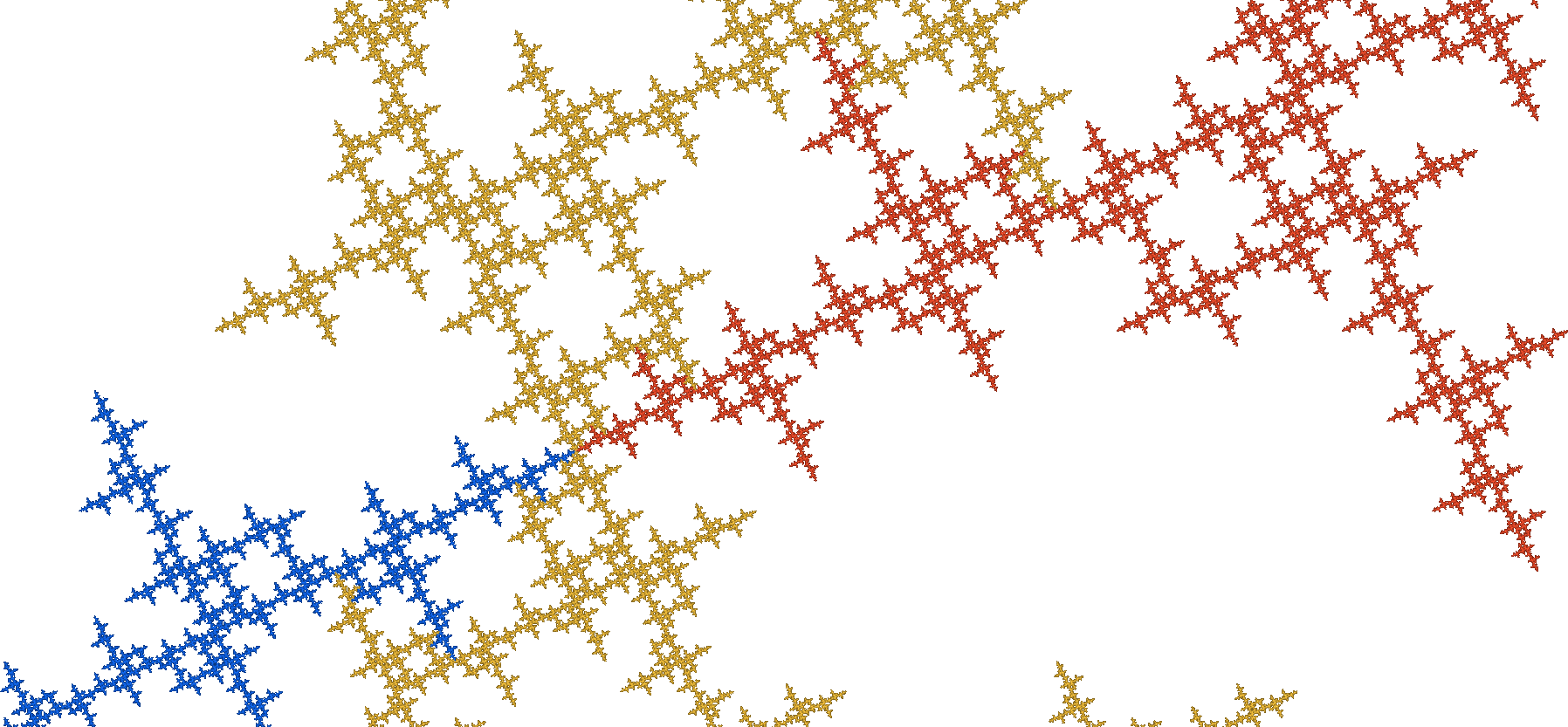}\vspace{4ex}\\
\includegraphics[width=0.9\textwidth]{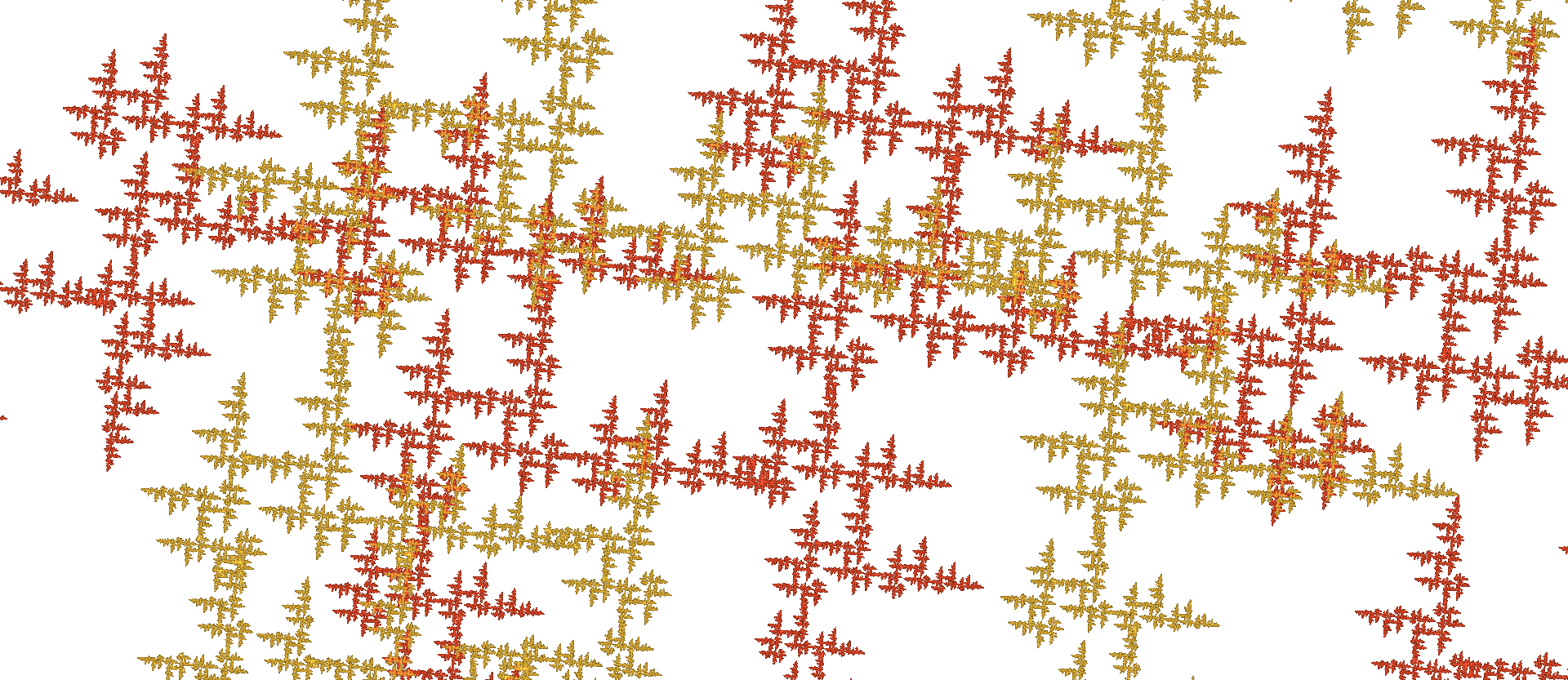}
\end{center}
\caption{Top: overlap region of the ``Crossings'' example in Figure \ref{examples}, generated by $f_0(z)=\frac{-iz}{2},\, f_1(z)=\frac{-1}{2}(z+1),\, f_2(z)=\frac{-1}{2}(z-1-i)$. Below: detail of ``Forest'', generated by $f_0(z)=\frac{i}{2}(z-5),\, f_1(z)=\frac{z}{2},$ and $f_2(z)=\frac{i}{2}(z-4+i).$ }
\label{crossing}
\end{figure}  

\paragraph{Classifying parameters.}
Since all fractals in this paper have the same fractal dimension, which parameters will distinguish them?
The Hausdorff dimension of the intersection of pieces can be determined from the leading eigenvalue of a matrix, as explained below. When $A$ is something like a polygon, then $A_j\cap A_k$ is something like an edge of $A.$ All spaces in Figure \ref{examples} have Cantor sets as edges, all with different dimensions, as Table \ref{paras} shows.  Except for ``Crossings'', this `boundary dimension' is larger than 1, while examples in Figure \ref{fsquares} can reach only values up to 1. Other parameters listed in Table \ref{paras} refer to the combinatorial structure of pieces and will be explained in Section \ref{neigh}.

\begin{table}[h!] 
\begin{center}
\begin{tabular}{|l|c|c|c|c|c|c|}\hline
name &Patches&Bubbles&Fireworks&Thicket&Crossings&Forest\\ \hline
topology&Cantor&Cantor&Cantor&conn.&conn.&conn.\\ \hline
holes&no&no&no&yes&yes&yes\\ \hline
line segments&no&no&no&no&yes&no\\ \hline
proper nbs.&90&84&36&89&7&58 \\ \hline
finite nbs.&15&0&17&0&2&7   \\ \hline
boundary dim&1.17&1.16&1.33&1.13&0.61&1.14  \\ \hline
max degree&19&22&21&16&7&13  \\ \hline
neighborhoods&3183&7521&954&6456&19&5938  \\ \hline
\end{tabular}\vspace{-2ex} 
\end{center}
\caption{Properties of the examples in Figure \ref{examples}, explained in the text.}  \vspace{-2ex}
\label{paras}\end{table}

\paragraph{``Crossings'' - our simplest example.} 
This is a typical connected fractal in our family, with many line segments. The interval $I$ generated by $f_1$ and $f_2$ is a kind of diagonal in Figure \ref{examples}, and $f_0(I)$ intersects $I,$ so that all pieces are connected with each other by lines, on every level (cf. \cite{bar}, chapter 8). While $A_1\cap A_2$ is a single point, magnification shows that the yellow part $A_0$ intersects the other two pieces in a Cantor set. Moreover, $A$ is not simply connected: it surrounds infinitely many holes of the same `quadratic' shape. Below we take this fractal to demonstrate our methods which for the other examples can be performed only by computer.
The two other connected examples in Figures \ref{examples}, \ref{thicket} and \ref{crossing} do not contain any line segments and surround holes of various shapes.

\section{The structure of the open set}\label{open}
In some cases, open sets $U$ fulfilling \eqref{osc} are easy to find.  In Figure \ref{fsquares}, a square can be taken as $U,$  and the interior of the convex hull of $A$ is another possible choice. Now we shall see that open sets can be  complicated and difficult to find. Since $U$ contains each $f_j(U),$ it contains all its smaller copies $f_{j_1}f{j_2}\cdots f_{j_n}(U),$ and this implies that the fractal $A$ is contained in the closure $\overline{U}$ of $U.$ 

\paragraph{Self-similar tiles.}
When the interior \, ${\rm int\,}A$ of $A$ is non-empty, then the property $A\subset\overline{U}$ implies that the open set is essentially unique, $U={\rm int\,}A$ (one can subtract a nowhere dense set which is invariant under the $f_k$ but this is not useful). Now there are many self-similar tiles $A$ with an intricate structure, see the examples in \cite{Me}. One of the simplest examples is a union of two squares: $A=[-1,0]\times [-1,0]\cup  [0,1]\times [0,1]$ which is a self-similar set with four pieces, with contraction factor $\frac12 .$ The open set  \, ${\rm int\,}A$ is disconnected. \smallskip

\paragraph{Disjoint pieces.}
We give a general approach and explain how a computer checks the OSC. This requires some concepts and notation. We start with the simplest case where the pieces $A_k=f_k(A)$ are all disjoint. They then have distance $\varepsilon>0$ from each other, for some $\varepsilon>0.$ The $\delta$-neighborhood of $A,$
\[ U=\{ z | \mbox{ there is a point $a$ in $A$ with } |z-a|<\delta \} \] 
with $\delta<\frac{\varepsilon}{2}$ is an open set which satisfies \eqref{osc}. For the proof we need only the contraction property of the $f_k,$ that is $|f_k(y)-f_k(z)|\le r_k\cdot|y-z|$ for all points $y,z$ with constant $r_k<1.$ This implies that each $f_k(U)$ is contained in the $\delta r_k$-neighborhood of $A_k.$ Thus these sets are disjoint and contained in $U.$\smallskip

\paragraph{Neighbors and edges.}
When the $A_k$ intersect, certain points of $A$ cannot belong to $U.$ If $a$ belongs to $f_j(A)\cap f_k(A),$ then $b=f_j^{-1}(a)$ cannot be in $U$ since this would imply that $f_j(U)\cap f_k(U)$ contains $a.$ This also holds for multiindices ${\bf j}=(j_1,j_2,...,j_n)$ and ${\bf k}=(k_1,k_2,...,k_n)$ which will be called \emph{words} or \emph{words of length $n$}. They denote compositions of contractions $f_{\bf j}=f_{j_1}\cdots f_{j_n}$ and small pieces $A_{\bf j}=f_{\bf j}(A)$ on the $n$-th level. The set 
\[  A\cap f_{\bf j}^{-1}f_{\bf k}(A) = f_{\bf j}^{-1}(A_{\bf j}\cap A_{\bf k}) \]
must be disjoint to $U$ for all words $\bf j,k$ with $j_1\not= k_1.$ (Otherwise $f_{\bf k}(U)$ contains a point of $A_{\bf j}\cap A_{\bf k}$ and thus intersects $f_{\bf j}(U),$ since $A_{\bf k}\subset \overline{f_{\bf k}(U)}.$ This would contradict \eqref{osc},)

It has become custom to call $f_{\bf j}^{-1}f_{\bf k}(A)$ a \emph{neighbor} of $A,$ and $f_{\bf j}^{-1}f_{\bf k}$ the corresponding \emph{neighbor map.} These are the sets which would be added if we extend the self-similar construction outwards, considering $A$ as a piece of a still larger `superset'. Since $A$ can be assumed to play the part of any piece $A_j$ or even $A_{\bf j},$ there can be lots of neighbors of $A.$ See for example \cite{BR,BM}, \cite{Lo,ST} for the important case of self-affine tiles, and \cite{DKV,SW} for the Levy curve. The intersection $f_{\bf j}^{-1}f_{\bf k}(A)\cap A$ of the fractal $A$ with a neighbor is considered an \emph{edge} of $A,$ a notation which is very natural for tilings. 

The closure of the union of all such edges was called the \emph{dynamical boundary} of $A$ by M. Moran \cite{Mo1}.  An open set $U$ must not contain any point of the dynamical boundary $\partial A$ of $A.$ If the dynamical boundary is equal to $A,$ then the open set condition fails. Moran conjectured that conversely, $\partial A\not= A$ implies the OSC. This problem is still unsolved. However, Bandt and Graf \cite{BG} proved that the OSC holds if and only if no sequence of neighbor maps $f_{\bf j}^{-1}f_{\bf k}$ can converge to the identity map.

\paragraph{Data from a discrete group.}
For our special class of IFS, the situation is simpler. From \eqref{maps} follows
\begin{equation}
f_j^{-1}(y)=s_j^{-1}(2y)-c_j\quad\mbox{ and }\quad   f_j^{-1}f_k(z)= s_j^{-1}s_k(z) +s_j^{-1}c_k -c_j \ .
\label{sym}\end{equation}
With induction it is easy to verify \cite{BG} that any neighbor map has the form
\[ h(z)=f_{\bf j}^{-1}f_{\bf k}(z) = s(z+c)\quad\mbox{ with } s\in S \mbox{ and integers $a,b$ such that } c=a+ib.  \]
Maps of this kind represent the translational and rotational symmetries of the integer lattice. They form a well-known crystallographic group. This is a discrete group, and the condition of \cite{BG} simply says that the OSC is fulfilled if and only if the identity is not a neighbor map. Moreover, since $A$ is bounded there can be only finitely many neighbor maps $h$ with $h(A)\cap A\not= \emptyset ,$ which are called \emph{proper neighbor maps.} \smallskip

\paragraph{What the computer does.}
This is how the IFStile package checks the OSC for a given IFS: all maps of the form \eqref{sym} which can be proper neighbor maps are determined. If the identity map is among them, OSC does not hold, otherwise OSC is satisfied. Thus the computer does not operate with any open set.  See Section 5 for details.

\paragraph{The central open set.}
Nevertheless, we want to know more about the structure of $U.$ The \emph{central open set} \cite{BR} is defined as
\begin{equation} U_c=\{ z|\, d(z,A)<d(z,h(A))\ \mbox{ for every neighbor map } h\} .
\label{cos}\end{equation}
Here $d(z,A)=\inf\{ |z-a|\, |\, a\in A\}$ is the distance of the point $z$ from the set $A.$ The central open set contains those points which are nearer to $A$ than to any neighbor of $A.$ Similar as for the case of disjoint $A_k$ above, one can verify \eqref{osc} for $U_c.$ Moreover, $U_C$ is non-empty if and only if $\partial A$ is different from $A.$ 
This definition makes it possible to construct natural open sets, with polygonal shape. In our new examples, infinitely many edges will be needed. 

\begin{Proposition}[A construction of open sets] \label{openset}\hfill\\
If $V\subseteq U_c$ is a non-empty open set, then the union $U=\bigcup f_{\bf \ell}(V)$ taken over all words $\bf \ell$ of arbitrary length, will be an open set fulfilling \eqref{osc}.
\end{Proposition}
 
{\it Proof. } The relation $f_j(U)\subset U$ follows from the definition of $U.$ Moreover $V\subseteq U_c$ implies that 
\begin{equation}\label{pf}
V\cap f_{\bf j}^{-1}f_{\bf k}(V)=\emptyset \mbox{ for every neighbor map }h=f_{\bf j}^{-1}f_{\bf k}
\end{equation} with words $\bf j, k$ of equal length with $j_1\not= k_1.$ This holds because the points of $h(U_c)$ are nearer to $h(A)$ than to $A.$ (Note that $h$ is an isometry and $h^{-1}$ is also a neighbor map.)  Now if $f_j(U)\cap f_k(U)$ would contain a point, by definition of $U$ the point would belong to $f_jf_{\bf \ell}(V)\cap f_kf_{\bf \ell '}(V)$
for some words $\bf \ell,\ell '.$ By enlarging one of the sets if necessary, we obtain words $\bf \ell$ and $\bf \ell '$ of equal length, which then contradicts \eqref{pf}. This completes the proof.
\hfill $\Box $\smallskip

\begin{figure}[h!] 
\begin{center}
\includegraphics[width=0.8\textwidth]{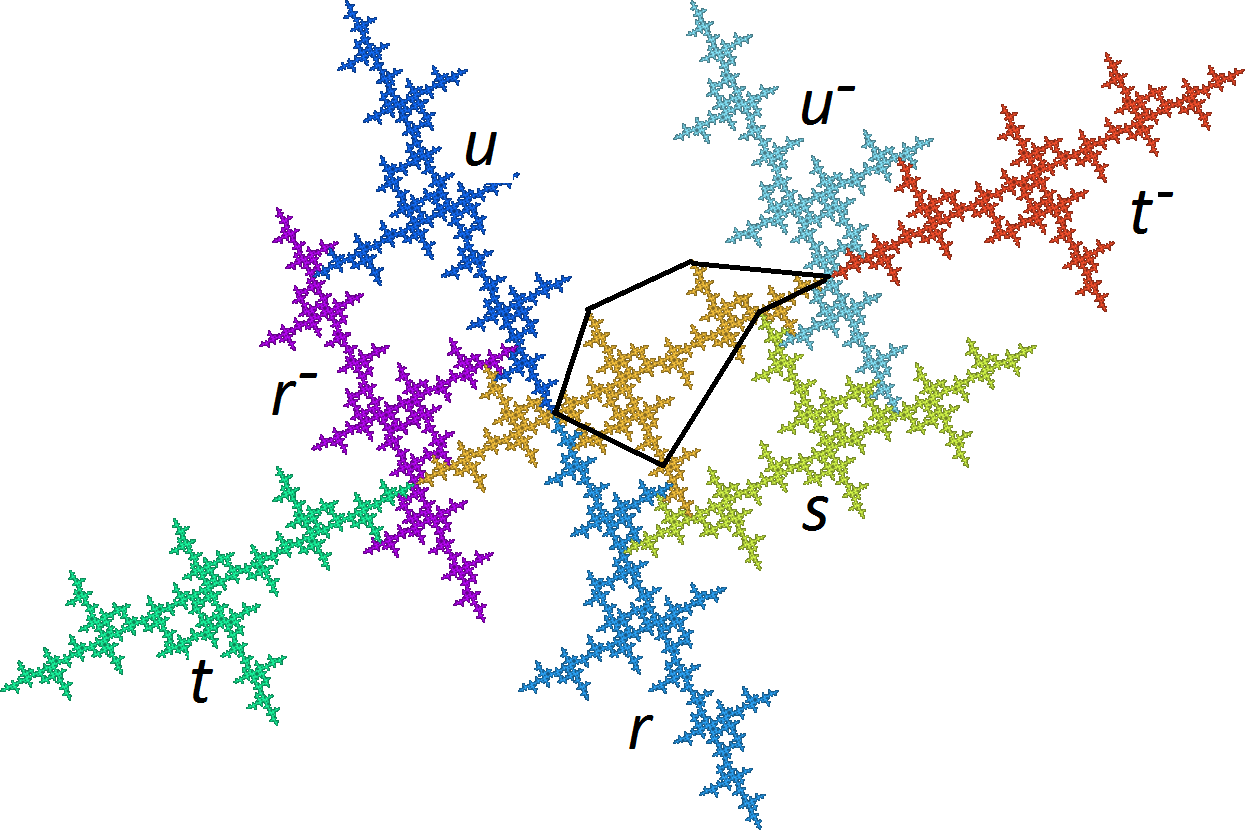} \vspace{-5ex}
\end{center}
\caption{``Crossings'' with its neighbor sets, and a polygonal set $V.$ Letter $t$ denotes $t(A).$}\label{crossingnb}
\end{figure}  

\paragraph{An open set for ``Crossings''.}
We apply the proposition to our simplest example.  Below we prove that ``Crossings'' has exactly seven neighbors, drawn in Figure \ref{crossingnb}. Obviously, the neighbors do not cover $A.$ So $U_c$ is not empty, but hard to calculate. As V we can take a small disk around some point of the central set $A$ which is far from all points of the neighbors. However, we would like to make $V$ as large as possible, so that the resulting set $U$ has a simple structure. A choice of $V$ as a fairly large polygon is indicated in Figure \ref{crossingnb}. The corresponding open set $U$ will be polygonal, with infinitely many sides and two connected components. $U$ is so simple since all neighbors do only touch $A$ without really crossing its line segments.

\begin{figure}[h!] 
\begin{center}
\includegraphics[width=0.6\textwidth]{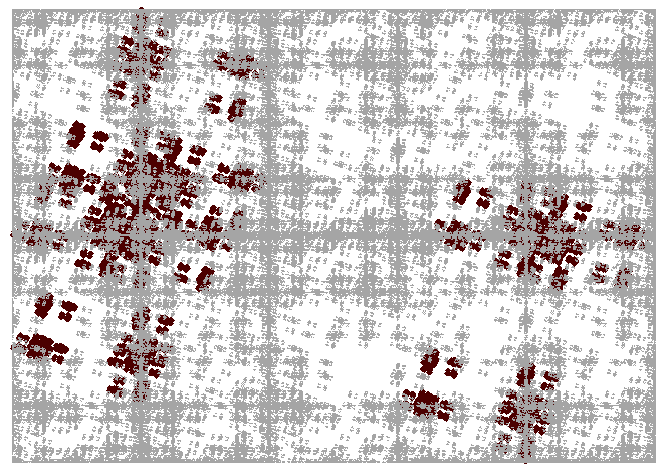}
\end{center}
\caption{``Fireworks'', as oversampled dark set, with its 36 proper neighbor sets drawn on top in light gray. A separating open set $U$ must be dense in the dark set while avoiding all neighbors. % Thus $U$ is extremely fragmented.
}\label{bubblenb}
\end{figure}  

\paragraph{A more complicated example.}
For comparison, we make the same experiment for the Cantor set ``Fireworks'' for which our program calculated 36 proper neighbor maps. The neighbors do not completely cover the set $A.$ However, there seem to exist only tiny connected sets $V.$ Any open set $U$ must be extremely fragmented. To study this self-similar set, we need details with observing window contained in an open set $U,$ thus with a huge magnification factor. This is why our detail pictures look so different from Figure \ref{examples}. For the remaining four examples, the situation is still worse.

\section{Parameters obtained from neighbor maps}\label{neigh}
\begin{figure}[h!] 
\begin{center}
\includegraphics[width=0.66\textwidth]{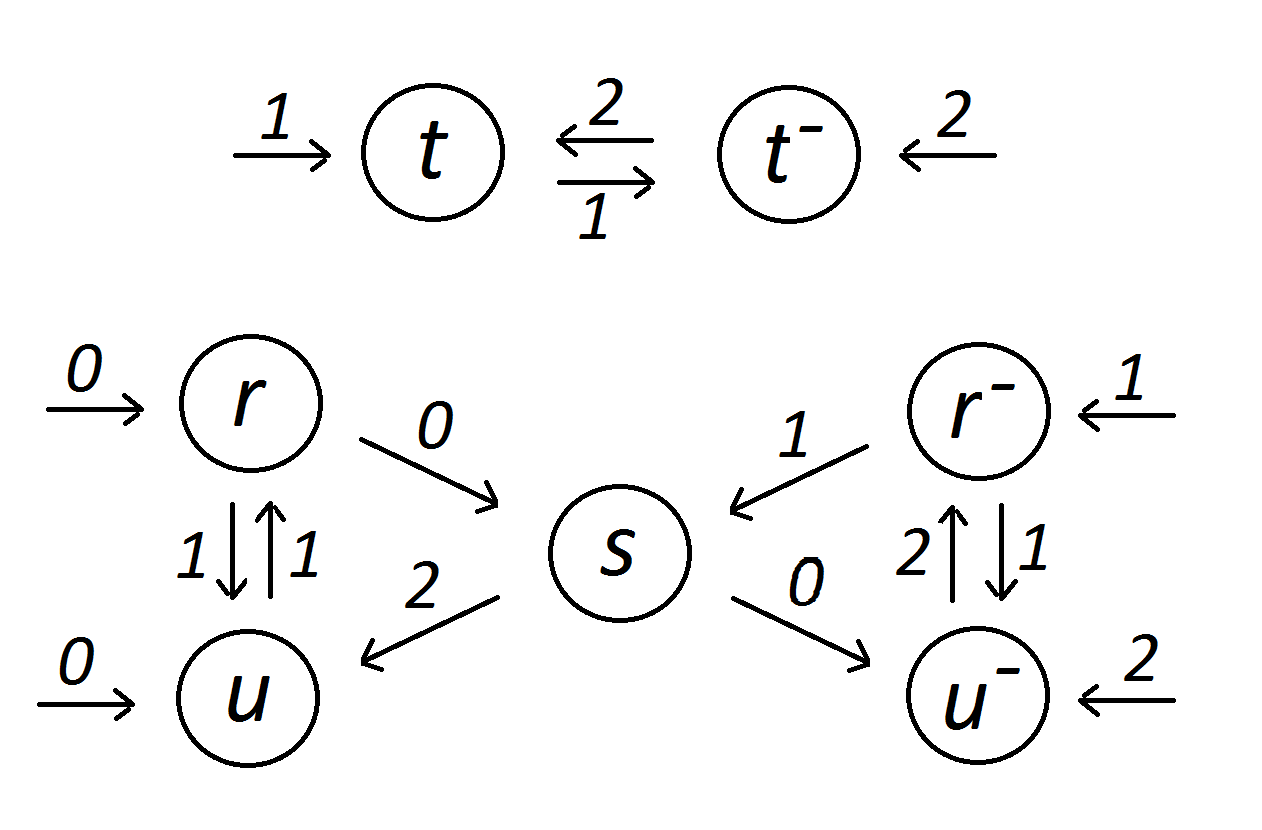}
\end{center}
\caption{The neighbor graph of ``Crossings'' has two parts describing finite and infinite neighbors.}\label{crossingraph}
\end{figure}

\paragraph{The neighbor graph.}
Neighbor maps do not only decide about the OSC. The proper neighbor maps form a graph which describes the topology of $A$ in the simplest possible way, and can be used to calculate various geometrical parameters for $A$ \cite{BM,BMT,DKV,Lo,ST,SW}. 
In the neighbor graph, every vertex is a neighbor map $h=f_{\bf j}^{-1}f_{\bf k}$ which represents the relative position of two intersecting pieces $A_{\bf j}, A_{\bf k}.$ For every pair of intersecting subpieces
$A_{{\bf j}j}, A_{{\bf k}k},$ there is 
 \[ \mbox{an edge labeled } (j,k) \mbox{ from neighbor map } h \mbox{ to } h'= f_j^{-1}hf_k= f_{{\bf j}j}^{-1} f_{{\bf k}k} .\]
Moreover, there are edges from $id$ to the initial neighbor maps $f_j^{-1}f_k$ with label $(j,k).$ Since $id$ is not a neighbor map, these initialization edges are drawn without initial vertex. In Figure \ref{crossingraph} we draw only the label $j$ since the second label $k$ can be found at the corresponding edge between inverse maps. Reflection at the vertical symmetry axis  associates each neighbor map with its inverse, and each label $j$ with the corresponding $k.$

Neighbor maps are often considered as types of intersecting pieces, and a self-similar set with finite neighbor graph is called \emph{finite type}. In our discrete setting, every $A$ is finite type. Only finitely many maps $h(z)=s(z)+v$ with $s\in S$ and integer translation vector $v$ can describe proper neighbors. Their number depends on the size of $A$ which in turn depends on the position of the fixed points of the $f_k.$ Thus it is possible to construct the neighbor graph by computer. Then it is checked whether $id$ is not a neighbor map, in which case the OSC holds true.

\paragraph{Neighbor maps for ``Crossings''.} In Figure \ref{examples} one can see that the pieces $A_1$ and $A_2$ of ``Crossings'' differ by a translation. Calculation shows $f_1^{-1}f_2(z)$ is the translation $t(z)=z-2-i$ which maps the central set in Figure \ref{crossingnb} to the neighboring set labeled with $t.$ Of course $f_2^{-1}f_1(z)$ is the inverse translation $t^-(z)=z+2+i.$ Each piece $A_{w1}$ for some word $w$ is mapped to $A_{w2}$ by a translation, which, standardized to the size and orientation of $A,$ is exactly $t.$ Now we consider the pieces $A_{12}$ and $A_{21}$ which contain the point $A_1\cap A_2$ on the second level. They also differ by a translation, but $A_{12}$ is mapped to $A_{21}$ by $t^-,$ because the $180^o$ rotation in $f_1,f_2$ reverses the direction. This fact, which the computer calculates as $f_{12}^{-1}f_{21}=t^-,$ is expressed by the upper row of Figure \ref{crossingraph}. 

 In Figure \ref{examples}, the piece $A_0$ is mapped into $A_1$ by a rotation around $-90^o.$ Calculating $f_0^{-1}f_1(z)$ we obtain $r(z)=-iz-i .$ The neighbor map between $A_1$ and $A_0$ is of course the inverse rotation $r^-(z)=iz-1.$ Similar calculations for pieces $A_0,A_2$ yield the rotations $u,u^-.$ On second level, there is one more neighbor map $s(z)=-z+1,$  a  $180^o$ rotation around the center of the biggest quadratic hole in Figure \ref{crossingnb}. It is self-inverse and corresponds to pieces $A_{00}, A_{11},$ or any pair of smaller pieces enclosing a `quadratic hole' in Figures \ref{examples} or \ref{crossing}. It turns out that there are no other neighbor maps because for each $h\in\{ t,t^-,r,r^-,u,u^-,s\} =\cal N$ and any pair $(j,k)$ of symbols in $\{ 0,1,2\},$ the map $h'=f_j^{-1}hf_k$ again belongs to $\cal N .$ We check that from each vertex, there is a path which leads to a cycle in the graph -- otherwise the caclculated neighbor would not intersect $A$ and had to be cancelled. The neighbor graph in Figure \ref{crossingraph} is complete. OSC holds because $id$ is not among the calculated neighbor maps.

As noted above, the algorithm must end after a finite number of steps. It was good luck, however, that our calculation ended on second level. The other examples require a computer.

\paragraph{Topology.} 
While the description of a self-similar set $A$ by its IFS is somewhat magic and intransparent, the
neighbor graph of a finite type fractal is a perfect mathematical tool. It gives an explicit description of the topology of $A$ as a quotient of the space of symbol sequences ${\bf j}=j_1j_2...$ The infinite paths of edges with labels $(j_1,k_1),(j_2,k_2),...$ exactly indicate those sequences $\bf j,k$ which are identified \cite{BM}. Of course, infinite paths in our finite graph can only exist when we have cycles. This is the reason for the check of neighbor maps mentioned above.

Let us study connectedness. The self-similar set $A=A_0\cup A_1\cup A_2$ is connected if one of the pieces intersects both of the other pieces, cf. \cite[chapter VII]{bar}. For ``Crossings'' $A_0\cap A_1$ is represented by $r$ and $r^-$, and $A_1\cap A_2$ by $t,t^-$ which already proves connectedness!  The upper part of Figure \ref{crossingraph} says that
$A_1\cap A_2$ is a point with two preperiodic addresses $1\overline{12}$ and $2\overline{21}.$ Standardized to $A,$ this means that the lower left endpoint with address $\overline{12}$ and the upper right endpoint with address $\overline{21}$ are boundary sets of $A.$ Here we could extend the construction outwards by a translate of $A,$ as is shown in Figure \ref{crossingnb}.

The lower part of Figure \ref{crossingraph} contains several cycles with a common vertex. This implies that $A_0\cap A_1$ and $A_0\cap A_2$ are Cantor sets. Two particular corresponding infinite paths are labelled $0\overline{1}$ on the left and $1\overline{12}$ on the right. Since the second address appeared already in the upper part, this shows that we have a point with three addresses starting with 0,1, and 2, that is, a common point of all three pieces, which can be seen in Figure \ref{examples}.  

\paragraph{Nonexistence of connected subsets.} 
If $A$ is not connected, it can still have large connected subsets. They are difficult to extract from the neighbor graph. One special case is easy: if there is an infinite path of edges $(j_1,k_1),(j_2,k_2),...$ where all labels are either 0 or 1, then the self-similar set $B$ with respect to $f_0,f_1$ is connected. This happens if both $f_0, f_1$ involve rotations of 0 or $180^o.$ 
Then $B$ is an interval, as in Figure \ref{fsquares}C.   For our data, this seems the only case to get connected subsets in a disconnected set $A.$ We verified the Cantor property of ``Patches'', ``Bubbles'', and ``Fireworks'' experimentally by mabnifying many pieces of the sets, and shall address the question in a subsequent paper.

\paragraph{Boundary dimension.}
Each neighbor map $h$ corresponds to an `edge' $H=h(A)\cap A$ of $A,$ that is, a piece of the dynamical boundary of $A.$ The neighbor graph directly determines a set of equations for the edges of $A.$ The edges themselves are self-similar sets which form a so-called graph-directed construction \cite{BM,BMT,DKV,Lo,ST,SW}. The OSC for $A$ implies the OSC for this construction. The Perron-Frobenius eigenvalue $\lambda$ of the adjacency matrix of the neighbor graph determines the Hausdorff dimension of the boundary pieces as $d=\log\lambda /\log r $ where $r$ is the common contraction factor of the maps $f_k,$ in our case $\frac12 .$ In this way, the dimension of the dynamical boundaries in Table \ref{paras} was determined by computer. 

We demonstrate the method for ``Crossings''.  Boundary equations can be read from the graph in Figure \ref{crossingraph}.  $R=r(A)\cap A$ denotes the boundary set generated by the map $r.$
\[ R= f_1(U)\cup f_0(S),\ U=f_1(R), \mbox{ and } S=f_2(U)\cup f_0(U^-)\ .\]
Note that $U^-$ is isometric to $U,$ namely $u(U^-)=u(u^-(A)\cap A) =A\cap u(A)=U .$ We substitute $U,U^-,$ and $S$ by images of $R$ and obtain
\[ R= f_1^2(R)\cup f_0f_2f_1(R) \cup f_0^2u^-f_1(R) \ .\]
Thus $R$ is a self-similar set with three similarity maps with factors $r_1=\frac{1}{4} , r_2=r_3=\frac{1}{8}.$ The Hausdorff dimension $\alpha$ of $R$ is determined by $\sum r_j^\alpha =1$ \cite{Mo,Fal,bar}. Putting $y=2^\alpha,$ we have  $y^3=y+2,$ with numerical solution $y\approx 1.5214.$ Thus $\alpha =\frac{\log y}{\log 2}\approx 0.6054$ is the dimension of $R,U,$ and $S.$ For Figure \ref{fsquares}B, a very similar calculation gives $3 (\frac{1}{4})^\alpha =1$ and  $\alpha =\frac12 \frac{\log 3}{\log 2},$ just half of the dimension of $A.$

For $T$ in Figure \ref{crossingraph}, as well as for any finite neighbor, the boundary dimension is zero. 
Different dimensions of infinite boundary sets can occur when the neighbor graph splits into different components. For Figure \ref{fsquares}H the set $A_0\cap A_1$ has dimension $\frac{2}{3}$ and $A_0\cap A_2$ has dimension $\frac{1}{3}.$ In our six examples, infinite boundary sets turned out to have the same dimension. When finite neighbors were omitted, the remaining neighbor graphs became irreducible, with exception of at most four ´transient' vertices in their initial part. 

\paragraph{Maximal degree and number of neighborhoods.}
Even when many neighbors exist, it is rare that all neighbors are realized at one piece $A_w.$ Actually, this happens for ``Crossings'' where $A_{011}$ has seven proper neighbors. The neighbor graph in Figure \ref{crossingraph} says that the piece $A_1,$ as well as all small pieces $A_{w1}$ with some 012-word $w,$ have neighbors corresponding to the maps $t$ and $r^-$ since these vertices are reached by an initial edge with label 1. All words $A_{w11}$ have three more neighbors corresponding to $t^-,s,$ and $u^-$ since these vertices are reached by a path labelled 11 and starting with an initial edge.
And, finally, $r,u$ can be reached by a path labelled 011.

A systematic counting of initial paths labeled with 012-words $w$ will also show us all other combinations of neighbors which are realized by pieces $A_w.$ As we have seen, the suffixes of $w$ must be included.  When a combination of neighbors  is realized by a piece $A_w$ which does not intersect the dynamical boundary, we call it a \emph{neighborhood.} 

Every neighborhood of ``Crossings'' has at least 2 neighbors since for each symbol 0,1,2 there are two initial edges labeled with this symbol. For instance, $\{ r,u\}$ is a neighborhood which is realized by  $A_{w20}$ for every word $w$ because the path 20 does not appear in the neighbor graph. Table \ref{paras} says that ``Crossings'' has 19 neighborhoods, ``Fireworks'' 954, and ``Bubbles'' even 7521. This number is small compared with $2^{\mbox{number of neighbors}},$ and seems a complexity parameter for $A.$

We can consider neighborhoods as vertices of a graph, in the same way as done for neighbors. We draw an edge with label $j$ from the neighbourhood of $A_w$ to the neighborhood of $A_{wj}.$ Since the $A_w$  do not meet the dynamical boundary, this graph turns out to be irreducible. It has a certain importance for the theory of microsets mentioned above.

The maximum degree in this neighbourhood graph coincides with the maximum number of neighbors which can occur at a piece $A_w.$ This is another complexity parameter for the set $A.$ It determines the largest possible density of the set, more precisely, of the canonical Hausdorff measure on $A.$ For ``Patches'', ``Bubbles'' and ``Fireworks'' we got almost equal degrees 19, 22, and 21, respectively.

Both maximum degree and number of neighborhoods describe the geometric network which arises from all $n$-th level pieces $A_w$ for large $n,$ when pieces are replaced points (their center of gravity, say), and points of intersecting pieces are connected.  There is no place for details here. Both parameters were determined by computer and listed in Tables \ref{fsq} and \ref{paras}.

\section{The interactive computer search}\label{last} 
In Section \ref{exa} we mentioned that the computer searches for new examples by a kind of random walk on the parameter space $\Pi$ of all IFS in our class. For this paper we used version 1.7.0.2 of IFSTile where a genetic algorithm is implemented. Random search alone will not lead to good results, however. The computer needs specific advice from the user, and the mathematical user should think a bit about the structure of the space $\Pi.$ One important issue is the existence of many IFS in $\Pi$ which generate essentially the same fractal $A.$

\paragraph{Equivalence of self-similar sets.}  
Two fractals $A,B$ are considered to be equivalent if there is a similarity mapping $g$ with $B=g(A).$ If $A$ is a self-similar set with respect to the IFS $\{ f_k| \, k=1,...,m\},$ fulfilling \eqref{selsi}, then $B$ will be the self-similar set with respect to the IFS $\{ gf_kg^{-1}| \, k=1,...,m\}$ because
\[ B=g(A)=g(\bigcup f_k(A)) = \bigcup gf_k(A) = \bigcup gf_kg^{-1}(B) \, .\]
In the complex plane $g$ can be a translation $g(z)=z+c,$ a rotation $g(z)=dz+c,$ or a orientation-reversing similarity map $g(z)=d\overline{z}+c$ with $d,c\in\CC .$ For example, reflections apply to the examples in Figure \ref{fsquares}, see \cite{bar,FC}. When $c,d$ have integer components, we remain in our class of IFS.

We can also consider affine maps $g(x)=M\cdot x+ v$ on $\RR^2,$ with a non-singular real $2\times 2$ matrix $M$ and $x,v\in \RR^2.$ Since the eye identifies a set with an affine image, we should consider $A$ and $B$ equivalent in this case, too. The problem is that the $gf_kg^{-1}$ need not be in our class of IFS, they may even fail to be a contracting maps.  However, there is one important special case. When $f_k(x)=rx+c$ with a real number $r$ then $gf_kg^{-1}$ has the same form.  In our class of IFS, this concerns examples with rotations only around 0 and 180 degrees. Figure \ref{fsquares}E, for instance, admits an affinely equivalent version in its fractal square family which is symmetric with respect to the axis $y=x.$ 

It should be noted that a \emph{permutation} of the maps of an IFS will of course not change the associated attractor $A.$ Thus equivalence also holds when the conjugacy $g$ is combined with a permutation of indices. This makes it more difficult to decide whether a long list of IFS already contains an equivalent of a new example.

\paragraph{Search options.} 
In our class of IFS  the Sierpinski gasket, Figure \ref{fsquares}E, and two similar examples have lots of affine equivalent versions. If we start with the Sierpinski gasket and try a random walk on $\Pi,$ then in $1/8$ of all IFS we select three rotations around 0 or 180 degrees, resulting in one of the four examples when the $c_k$ are not on a line. Since these examples fulfil the OSC while other data do not, we can have 90 percent of our results equivalent to one of the four examples.  This is a typical difficulty in our random search.

Below we show how to filter out equal examples. A perfect solution has not yet been found. One has also to avoid to visit too often exactly the same IFS. Of course one can keep a list of visited places in $\Pi$ and compare each new IFS with this list, but this can become expensive when we visit millions of places. Since many IFS do not satisfy the OSC, it seems better to save only the list of successful examples.

In our case, we avoided getting the four trivial examples by fixing the angle 90 degrees for the first mapping. It would also be possible to define the first translation vector as $c_0=0,$ thus reducing the dimension of the parameter space by two (real and imaginary part of $c_0$). However, with fixed $c_0$ we got less results than with variable $c_0.$ Apparently small $c_0$ and OSC force $c_1,c_2$ to be large. 

On the whole, the search was done with relatively small translation vectors, $a_k,b_k$ mostly between -10 and 10. Larger translations are usually associated with large neighbor graphs requiring a lot of computation. The search has to be done in such a way that any single complicated example is handled within milliseconds. If the IFS is too complex, it must be abandoned.

\paragraph{Complexity parameters.} 
For our class of fractals, all neighbor maps are in the rotational symmetry group of the integer lattice, and all calculations can be done with integers. All results are rigorous. No numerical approximations are used. When the program says `OSC is true' then this is correct!  However, when the program abandons a data set, it may still fulfil the OSC, with a large number of neighbors. 

It is possible to prescribe a complexity bound for the random search. Any IFS which exceeds this bound will be skipped.
The complexity parameter is not the exact number of proper neighbors since this would take too much time. Instead, a bound for the size of $A$ is estimated from the data of the IFS, and this gives an estimate for the number of possible neighbors. We shall not go into details here.

\paragraph{Sort and select results.}
When the program is stopped, we have a list of IFS which fulfil the OSC. Every single fractal, or matrices of $m\times n$ fractals, can be drawn on the screen and studied by the user. However, it may happen that there are too many examples in the list to be studied individually. And many of the examples may be not interesting. Thus a selection has to be made. A preselection is performed by cancelling examples where all pieces of $A$ are disjoint -- although the OSC holds in this case, as shown in Section \ref{open}.

For an efficient selection, the program can determine a number of parameters for each example, as indicated in Section \ref{neigh}. This includes Hausdorff dimension of $A,$ boundary dimension of $A,$ number of possible neighbors, number of edges in the neighbor graph, complexity estimate, topological properties like connectedness and also certain moments of the self-similar measure. One can sort according to these parameters, to see for instance examples with large boundary dimension. 

It is also possible to filter with respect to certain parameters, and cancel IFS which seem not interesting.
Among 100000 IFS in our list, more than 90000 were like Figure \ref{fsquares}C, disconnected with intervals. After omitting cases with exactly two neighbors, less than 10000 examples were left for study.
Another option is to allow only a certain number $n$ of examples with the same value of boundary dimension, or of moments.  This can considerably shorten the list of results by avoiding repetition of equivalent cases, but some interesting examples can get lost.

Thus, although the IFS tile finder is a very nice program, it will not run automatically. It must be used in interactive mode, with repeated experiments, and the mathematician's wisdom must combine with the computing power of the machine.

\section{Conclusion.} 
Self-similar sets are often considered as a class of toy models, with clear structure, easy to understand. Here we have shown that even for very simple similarity maps, there is a great variety of self-similar sets with unexpected complexity.  Interactive computer search is the appropriate method to find new examples, study their properties and classify them. On the mathematical side, the structure of separating open sets and the topological and metric properties of the fractals were studied with the help of neighbor maps.

\begin{acknowledgements} This research was supported by the German Research Council (DFG), project Ba 1332/11-1. When the paper was designed, the first author enjoyed the hospitality of Institut Mittag-Leffler, Stockholm, Sweden, within the program `Fractal Geometry and Dynamics'.
\end{acknowledgements}

\end{document}